\documentclass[11pt,leqno]{article}
\usepackage{amsthm,amsfonts,amssymb,epsfig,graphics,amsmath,eufrak,oldgerm}
 \usepackage[latin1]{inputenc}\relax

\newcommand{\RR}{{\mathbb R}}
\newcommand{\NN}{{\mathbb N}}
\newcommand{\ZZ}{{\mathbb Z}}
\newcommand{\TT}{{\mathbb T}}

\newcommand{\EE }{{\mathbb E}}

\newcommand{\FF}{{\mathbb F}}

\newcommand\cN{{\cal  N}}

\newcommand\cT{{\mathcal T}}

\newcommand{\ls}{{\ \lesssim \ }}

\def\D{\partial }
\newcommand\adots{\mathinner{\mkern2mu\raise1pt\hbox{.}
\mkern3mu\raise4pt\hbox{.}\mkern1mu\raise7pt\hbox{.}}}

\newcommand{\im}{{\rm Im }\, }
\newcommand{\re}{{\rm Re }\, }

\newcommand{\mez}{\frac{1}{2}} 
\newcommand{\qt}{\frac{1}{4}} 
\newcommand{\tqt}{\frac{3}{4}} 
\newcommand{\uE}{\underline{E}}

\newtheorem{theo}{Theorem}[section]
\newtheorem{prop}[theo]{Proposition}
\newtheorem{cor}[theo]{Corollary}
\newtheorem{lem}[theo]{Lemma}

\newtheorem{rem}[theo]{Remark}

\numberwithin{equation}{section}

\title{Instabilities in  Zakharov Equations }

 \author{Thierry Colin  
 and Guy M\'etivier \footnote{MAB,  Universit\'e de Bordeaux I,
33405 Talence cedex, France. Email : 
Thierry.Colin@math.u-bordeaux1.fr, \  Guy.Metivier@math.u-bordeaux1.fr} }

\begin{document}

\maketitle

\begin{abstract}

In \cite{LPS}, F.Linares, G.Ponce, J-C.Saut have proved that  a non-fully dispersive  Zakharov system arising in the study of Laser-plasma interaction, is locally well posed in the whole space, 
for fields vanishing at infinity. Here we show that in the periodic case, seen as a model for 
fields non-vanishing at infinity, the system develops strong instabilities 
of Hadamard's type, implying  that the Cauchy problem is strongly ill-posed.

\end{abstract} 
\tableofcontents

\section{Introduction}

The construction of  powerful lasers allows new
experiments where hot plasma are created in which laser beams can
propagate. The main goal is to simulate in a laboratory 
nuclear fusion by inertial confinement. We  need some precise and reliable models of laser-plasma 
interactions in order to produce numerical simulations than are useful in order to understand the
experiments. The kinetic-type models are the more precise ones but the
cost in term of computations is exorbitant and no physically relevant situation
for nuclear fusion can be simulated in this context. Therefore, we use a
bi-fluid model for the plasma : we couple two compressible Euler systems
with Maxwell equations. Even under this form, it is not possible to 
perform direct computations because of high frequency motions and the 
small wavelength involved in the problem. At the beginning of the 70's,
Zakharov and its collaborators introduced the so-called Zakharov's
equation \cite{zakharov} in order to describe the electronic plasma
waves. Basically, the slowly varying envelope of the electric 
field $E=\nabla \psi$ is coupled to the low-frequency variation of
the density of the ions $\delta n$ by the following equations which are written
in a dimensionless form :
\begin{eqnarray}\label{1.1}
\left\{\begin{array}{l}
i\partial_t \nabla \psi +\Delta (\nabla \psi)=\nabla \Delta^{-1}
\text{div}(\delta n\nabla \psi),\\\\
\partial_t^2 \delta n-\Delta \delta n =\Delta(|\nabla \psi|^2).
\end{array}\right.
\end{eqnarray}
Of course, variations of this systems exists (see \cite{sulem1} for example).
For laser propagation, one uses the paraxial 
approximation and the Zakharov system reads

\begin{equation}
\label{1}
\left\{ \begin{aligned}
& i( \D_t  + \D_z)E + \Delta_x E = n E,\\
& (\D_t^2 - \Delta_x) n = \Delta_x \vert E \vert^2
\end{aligned}
\right.
\end{equation}
The space variables are  $(z, x)$, 
$z \in \RR$ and $x \in \RR^2$; $z$ is the direction of propagation of the laser beam
 and $x$ are the directions transversal  to the propagation. See \cite{Z2}
or \cite{sulem1} for a symmetric use of this kind of model
for numerical simulation.

We consider the Cauchy problem for \eqref{1} with initial data
\begin{equation}
\label{1bc}
\left\{ \begin{aligned}
   E_{| t = 0}  =   E_0, & \\
   n_{| t = 0}  =   n_0  , &  \qquad  \D_t  n_{| t = 0}  =   n_1. 
\end{aligned}
\right.
\end{equation}
The   existence theorem (see \cite{sulem1, GTV, OT} and references therein)
for the classical Zakharov system, that is when $\Delta_x$ is replaced by $\Delta_{(z, x)}$, 
does not apply.  
In \cite{LPS}, it is proved that the Cauchy problem for 
\eqref{1} is  well posed, locally in time, for data in suitable  Sobolev spaces. 
The proof is based on dispersion estimates. For periodic data, these dispersion estimates 
are not valid. This is a well known phenomena, even in the simpler case  of  Schr\"odinger 
equations. However, the new phenomena  here is that the consequences of this  lack of dispersive effects  are much more dramatic since it implies strong instabilities of Hadamard's type, 
so that the Cauchy problem for periodic data is strongly ill-posed in Sobolev spaces. 

For numerical experiments, one uses mainly periodic boundary conditions: one considers that the
plasma is infinite and has a periodic structure. For this kind of
application, it is quite reasonable to consider that
$E \sim \underline E \ne 0  $ at infinity. Our result has therefore a practical application and means that the paraxial
approximation is not a good model in this case: one should add the longitudinal dispersion.

\bigbreak 

We look for solutions $U = (E, n)$ of \eqref{1}, which are periodic in $x$, with period $2 \pi$
in $x$
and periodic in $z$ with period $2 \pi   Z$, where $Z$ is arbitrary. We denote by $\TT$ the 
corresponding torus  $\RR/ 2 \pi Z  \times (\RR/ 2 \pi)^2$. 

We consider the constant solution 
\begin{equation}
\underline U =  (  \underline E,  0) , \quad   \underline E \ne 0, 
\end{equation} 
which of course does not belong to the spaces used in \cite{LPS}, and 
we prove that this solution is strongly unstable.

\begin{theo}
\label{mainth}
For all $s$, there are families of  solutions  $U_k  = \underline U + 
(e_k, n_k )$, in $C^1([0, T_k] ; H^s(\TT))$
such that 
\begin{eqnarray}
&& \|  e_k(0) , \,  n_k(0) , \, \D_t n_k (0) \|_{H^s(\TT)}  \ \to \ 0,  
\\
&& T_k \to 0, 
\\
&&     \|  e_k( T_k) , \,  n_k(T_k)  \|_{L^2(\TT)}  \ \to \  \infty . 
\end{eqnarray}

\end{theo}

This nonlinear instability result is pretty strong:  not only the 
\emph{amplification}    $ \| u(T) \|_0  / \| u(0) \|_s  $ is arbitrarily large, in 
arbitrarily small time $T$, with arbitrary loss of derivatives $ s$, but there is 
an effective \emph{blow up} of the $L^2$ norm.


\section{The instability mechanism}

Our  construction is based on  the  analysis of the \emph{dispersion relation} for the
Zakharov system. 
 Consider the linearized equations around $(\underline E,   0 )$ : 
\begin{equation}
\label{2}
\left\{ \begin{aligned}
  i( \D_t + \D_z)  e + \Delta_x   e  -   \underline E  \, n  & =  f  ,\\
  (\D_t^2 - \Delta_x)   n -  \Delta_x ( 2 \re  \underline  E \  \overline e)  & =  g 
\end{aligned}
\right.
\end{equation}
With $(e, \overline e ,   n)$ as unknowns  the  system reads: 
\begin{equation}
\label{2b}
\left\{ \begin{aligned}
  - i( \D_t + \D_z)  e -  \Delta_x    e  +   \underline E  \, n  & =  -  f  ,\\
  i( \D_t + \D_z)  \overline e -  \Delta   \overline e  +    \underline E  \, n  & =  - \overline f  ,\\
  (\D_t^2 - \Delta_x)   n -         \underline  E \  \Delta_x \overline e  - 
    \overline{ \underline  E}  \   \Delta_x  e  & =   g . 
\end{aligned}
\right.
\end{equation}
Denoting by $(\tau, \zeta, \xi)$ the frequency variables dual to $(t, z, x)$, its symbol is  
\begin{equation}
\label{3}
 \begin{pmatrix}
  ( \tau+ \zeta) + \vert \xi \vert^2    & 0 & E   \\
     0  &  - ( \tau+ \zeta) + \vert \xi \vert^2 & \overline E 
     \\
     \vert \xi \vert^2 \overline E &  \vert \xi \vert^2   E  & \vert \xi \vert^2 - \tau^2 
\end{pmatrix}
\end{equation}
and the relation dispersion is $P = 0$, where $P$ is the determinant of the system, that is  
\begin{equation}
\label{4}
 P = ( \vert \xi \vert^2 - \tau^2 ) \big( \vert \xi \vert^4  - (\tau+ \zeta)^2\big) - 2 
 \vert \underline E \vert^2 \vert \xi \vert^4  = P_0 - 2 
 \vert \underline E \vert^2 \vert \xi \vert^4. 
\end{equation}
The remark is that  for $(\zeta, \xi)$ real,  $P_0$ has 4 real roots in $\tau$ 
\begin{equation}
\label{5}
- | \xi|,  \quad +| \xi | ,  \quad - \zeta - \vert \xi \vert^2 , \quad - \zeta + \vert \xi \vert^2,  
\end{equation}
 with an intermediate  double root  when $0 < | \xi  | = -\zeta - | \xi|^2 $. 
 Note that $P_0$ is   of degree 6 in $\xi$ while the perturbation 
 $-  \vert \underline E \vert^2 \vert \xi \vert^4 $ is of degree $4$ and negative. 
 Therefore, for $\xi$ large and $ \zeta  =  - | \xi|  - | \xi| ^2  $, the double root of 
 $P_0$ is perturbed in two conjugated \emph{complex roots}. 
 More precisely, for 
\begin{equation*}
\label{6}
| \xi |   \gg 0, \quad   \zeta  =   - | \xi|  - | \xi| ^2  \quad \mathrm{and} 
\quad \tau =  | \xi | (1 + \sigma ) , 
\end{equation*} 
the determinant $P$ is 
\begin{equation}
\label{8}
P = -  |\xi |^5  \Big(  \sigma^2 ( 2   - \sigma/ | \xi | ) ( 2   +  \sigma) + 2 
 \vert \underline E \vert^2 /  \vert \xi \vert  \Big)  . 
\end{equation}
The implicit function theorem  shows that there are two non-real roots 
\begin{equation}
\label{10}
\tau = \xi  \pm i  \frac{| \underline E | }{\sqrt 2}  | \xi |^{\mez}  + 0 (1) . 
\end{equation}

This means  that waves at frequency $(\zeta, \xi) $ 
with $\zeta  =   - | \xi|  - | \xi| ^2  $ are amplified by the exponential factor 
\begin{equation}
\label{exp}
e^{ \gamma t | \xi |^\mez} ,   \quad \gamma = \frac{| \underline E | }{\sqrt 2}  > 0 . 
\end{equation}
This implies that the Cauchy problem for the linearized equations \eqref{2}  is ill-posed 
in $H^\infty$ : there are Cauchy data in $H^\infty$ such that the homogeneous problem 
with $f = g =0$  has no solution 
in $C^0([0, T]; H^{- \infty})$. 

The goal of this paper is to translate this \emph{spectral instability} into a nonlinear instability 
result for the Zakharov system \eqref{1}.

\begin{rem}
\textup{How is it that this spectral instability does not intervene in the analysis of \cite{LPS}? 
The first answer is that the condition $ \underline E  \ne 0$ is crucial for $\gamma$ to be positive.
In their case, where solutions vanish at infinity, linearizing the equation around non-vanishing constants has no real significance.  However, the symbolic calculus above also makes sense
in the case of variable coefficients and one   expects  that the dispersion 
relation $P = 0$, with $\underline E$ replaced by $E(t, z, x)$, which still has non-real roots, 
should play an important role in the analysis.  For instance, the symbolic analysis 
appears when one replaces the plane wave analysis used for constant coefficients, 
by  \emph{geometric optics} expansions associated to localized 
wave packets. In this case, for a wave packet with mean frequency  $(  - | \xi|  - | \xi| ^2, \xi) $ 
an exponential amplification 
similar to \eqref{exp} is   expected. But the group velocity in $x$ of this packet is of order 
$ 2 \xi $; therefore if $E$ is confined (think of it as compactly supported) the time of 
amplification is short (typically $O(| \xi |^{-1})$ ) so that the overall effect of the 
amplification is bounded. Of course, this is just a very rough explanation, but it is 
rather intuitive. The detailed   balance between amplification and localization 
is indeed given by the dispersive estimates proved in \cite{LPS}.  }

\end{rem}

\begin{rem}
\textup{The system can be reduced to first  order in $t$, introducing $(\D_x e, \D_t n, \D_x n) $ as  unknowns, but it is not first order in $x$, because of the Schr\"odinger part of the system. 
However,  there is a good analogy with the analysis of 
weakly hyperbolic system. 
Indeed, the analysis of the symbol \eqref{3} shows that  when  for $\zeta = -  | \xi |^2 - | \xi | $,
there is a double eigenvalue with a $2 \times 2 $ Jordan block.   
The exsitence of  non-real  eigenvalues  \eqref{10},   simply means that the natural analogue
of the Levi condition for first order system is not satisfied. Pursuing the analogy, 
the exponential growth \eqref{28} indicates that the  Cauchy problem should be well posed in Gevrey classes $G^s$ for $s \le 2$. }

\end{rem}


\section{Scheme of the proof}

It is certainly sufficient to prove the theorem with  functions 
of $ x = (x_1, x_2)$ independent of $x_2$. To simplify notations,   
we assume from now on that $x $ is one real variable.  
Consider spatially periodic solutions of \eqref{1}, with period 
$2 \pi$ in $x$ and $2 \pi   Z$ in $z$. 
Moreover, 
 we look for solutions  $n$ and $E$ of the form  
\begin{equation}
\label{21}
\begin{aligned}
  n  &=    n (k x -  m z , t ) \\
  E &=  \underline E   +    e (k x -  m z , t )
\end{aligned}
\end{equation}
with new functions $n(\theta, t) $ and $e(\theta, t) $ $2\pi$ periodic in $\theta$. 
For the functions to be $2 \pi$ periodic in $x$ and $2 \pi / Z$ periodic in $z$, 
it is sufficient that
\begin{equation}
\label{22}
k \in \NN, \quad   m  Z    \in    \NN. 
\end{equation}
To be close to the  unstable frequencies, we require that 
$\vert m - k - k^2 \vert \ll \sqrt k $ and therefore we choose 
$m \in      \NN / Z   $ such that 
\begin{equation}
\label{23}
 \quad (k^2 + k) -   1/Z   <     m   \le  (k^2 + k)  . 
\end{equation}

The new equations read 
\begin{equation}
\label{24}
\left\{ \begin{aligned}
& i( \D_t -  m \D_\theta ) e +  k^2 \D^2_\theta  e  -  \uE  n  =   n e  ,\\
& (\D_t^2 - k^2 \D^2_\theta )   n -   k^2 \D^2_\theta (  \overline {\uE} \, e\,+ \,  \uE\,  \overline e ) 
=   k^2 \D_\theta^2 
\vert e \vert^2 .
\end{aligned}
\right.
\end{equation}
With $U := {}^t (e, n)$,  write it as
 \begin{equation}
\label{25}
L_k( \D_t, \D_\theta) U =    N_k (  U) 
\end{equation}
where $L_k$ is  the linear operator defined in the left hand side of \eqref{24}, 
and $N_k(u)$ the quadratic term in the right hand side.  

\bigbreak

The first step concerns the homogeneous equation 
\begin{equation}
\label{linhom}
L_k   U = 0 , 
\end{equation}
which  is studied using Fourier series expansions in $\theta$. 
The choice \eqref{23}
and the spectral analysis of Section 2 and the choice \eqref{23} 
 imply  that for $k$ large, the harmonic $1$ is unstable :  
 
\begin{prop}
\label{p1}  There is $k_0$ such that for $k \ge k_0$, 
there are solutions 
 $U^a = (e^a , n^a)  $  of $\eqref{linhom}$  such that 
   \begin{equation}
   \label{37} 
\left\{   \begin{aligned}
 e^a  & =    \hat e^a_{1}(t)    e^{ i    \theta }   + \hat e^a_{-1}(t) e^{ - i \theta} 
\\  n^a & =  \sinh \big(   t \sigma \big)  \cos \big(  t \re \lambda  + \theta \big)  
 \end{aligned} \right. 
 \end{equation}
with 
\begin{equation}
\label{38a}
\hat e^a_{\pm 1}  (t)  =    \big( e^a_{\pm 1  , + }   e^{  t \gamma } +    e^a_{\pm 1 ,- }  
e^{ - t \gamma } \big)  e^{ i    t  \lambda}   ,   
\end{equation}
where the parameters $\lambda , \sigma  , e_{\pm 1, \pm} $ depend on $k$, 
$\lambda$ and $\sigma $ being  real positive  and, 
  as $k \to + \infty$, there holds : 
\begin{equation}
\label{39a}
e^a_{+ 1 , +}  \sim  -  i \overline \uE  / 4 \sigma , \quad  
e^a_{+ 1 , -}  \sim  -  i  \uE  / 4 \sigma , \quad 
e^a_{ - 1, \pm } = O(k^{-2})   .   
\end{equation}
\begin{equation}
\label{310a}
 \lambda \sim k, \qquad \sigma \sim | \uE |  \sqrt { k^ / 2}. 
\end{equation}
\end{prop}
The proof is given in Section~4. 

\medbreak
 Next, we  consider $\delta U^a$ as a first approximation of the 
 solution of \eqref{25} to construct,  with $\delta$ a small parameter to be chosen. 
 More precisely  look for solutions of \eqref{25} as 
\begin{equation}
\label{27y}
U = \delta ( U^a +   u  ) , \quad  u =  (e, n)
\end{equation}
with the same initial data as $\delta U^a$. 
Because the nonlinearity is exactly quadratic, the equation for $u$ reads
 \begin{equation}
\label{28y}
L_k( \D_t, \D_\theta) u  =    \delta N_k (  U^a + u ) \,, \quad  e_{\vert t = 0} =
 n_{\vert t = 0}  =  \D_t n _{\vert t = 0}  = 0 . 
\end{equation}
This equation is solved by Picard's iteration  and  the main 
step is to solve the linear equation 
\begin{equation}
\label{26}
L_k  U =  F  ,    \quad  e_{\vert t = 0} =
 n_{\vert t = 0}  =  \D_t n _{\vert t = 0}  = 0 . 
\end{equation} 
in   Banach spaces  which are  also well adapted  to the nonlinearity. 
The choice of these spaces, more precisely of their norm, is technical
and dictated by the computations detailed in the next sections. 
We just give here their definition. 

For a periodic function $v$ of $\theta$, we denote by 
$\hat v_p$ its Fourier coefficients so that 
\begin{equation}
\label{314}
v = \sum_{p \in \ZZ}  \hat v_p  e^{ i p \theta} . 
\end{equation}
The  first Fourier coefficient    $\hat e_1$    plays a special role and we use the notations
\begin{equation}
\label{42}
e(t, \theta)  = \hat e_1 (t) e^{ i \theta}  + e' (t, \theta). 
\end{equation}
For  
$ s \ge 1$ and $T > 0$ ,   se denote by 
$\EE^1(T)$ the space of  $u = (e, n)$  with $n$ real valued,  such that 
\begin{equation}
\label{316} 
 e \in C^0([0, T] ; H^{s+2} ) \cap C^1([0, T] ; H^s) , \qquad 
 n \in C^1([0, T] ; H^s) 
\end{equation}
 equipped with the norm
 \begin{equation}
 \label{317} 
 \begin{aligned}
 \| u \|_{\EE^1(T)} =   \sup_{ t \in [0, T]}  \ & e^{ - \sigma t } \Big\{   
 k^\mez   \vert \hat e_1  (t) \vert  +  k^{- \mez}      \vert \D_t \hat e_1  (t) \vert   
 +  k^{\tqt}   \Vert e'  (t) \Vert_{H^{s+2} } 
 \\
 &+    
k^{-\mez}  \Vert \D_t  e' (t) \Vert_{H^s}  + \Vert n (t) \Vert_{H^s} +  k^{-1}    \Vert \D_t n (t) \Vert_{H^s} 
 \Big\}
 \end{aligned}
 \end{equation}
 where $\sigma$ is defined at Proposition~\ref{p1}. The norm depends on $k \ge 1$ and $s$, 
 but, to lighten  the text,  we do not mention this dependence explicitly in the notations. 
 
 We denote by $\EE^2(T)$ the same space \eqref{316}, equipped with the norm
  \begin{equation}
 \label{318} 
 \begin{aligned}
 \| u \|_{\EE^2(T)} =     \sup_{ t \in [0, T]} & \   e^{ -  2 \sigma t } \Big\{   
 k    \vert \hat e_1  (t) \vert  +         \vert \D_t \hat e_1  (t) \vert   
 +  k    \Vert e'  (t) \Vert_{H^{s+2} } 
 \\
  +    & 
k^{-\qt}  \Vert \D_t  e' (t) \Vert_{H^s}  + 
k^\mez \Vert n (t) \Vert_{H^s} +  k^{-\mez}    \Vert \D_t n (t) \Vert_{H^s} 
 \Big\}. 
 \end{aligned}
 \end{equation}
 There are two differences between \eqref{317} and \eqref{318} : first the weight 
 $e^{ - \sigma t } $ is replaced by $e^{ - 2 \sigma t } $ and second all the powers of  
   $k$ in the coefficients are increased, at least by a factor $\qt$. In particular, 
   \begin{equation}
   \label{compE}
    \| u \|_{\EE^1(T)}   \le     k^{-\qt}  e^{ \sigma T}  \| u \|_{\EE^2(T)} 
   \end{equation}

For the right hand sides, we denote by $\FF^2(T)$ the space of $F = (f, g)$ with $g$ real valued such that 
\begin{equation}
\label{319} 
 f \in   C^1([0, T] ; H^s) , \qquad 
 g \in C^0([0, T] ; H^s)\quad  \mathrm{with} \ \hat g_0 = 0 , 
\end{equation}
 equipped with the norm
 \begin{equation}
 \label{320} 
 \begin{aligned}
 \| F \|_{\EE^1(T)} =   \sup_{ t \in [0, T]}  \  e^{ - 2 \sigma t } \Big\{   
   k^{\mez }   \Vert f  (t) \Vert_{H^{s} } 
&  +    
k^{-\mez}  \Vert \D_t  f (t) \Vert_{H^s} 
\\
& +  k^{-\tqt }    \Vert g (t) \Vert_{H^s} 
 \Big\}
 \end{aligned}
 \end{equation}

The   next three results justify the choices of these norms. 
We assume that the parameter $s \ge 1$ is fixed. 

The first estimate is an immediate consequence of  Proposition~\ref{p1} and \eqref{39a}
\eqref{310a}. 
\begin{lem}
\label{l2}
There is a constant $K^a$ such that for all $k \ge k_0$ and all $T \le 1$, 
the approximate solution $U^a$  of Proposition~$\ref{p1}$ satisfies 
\begin{equation}
\label{321}
\|  U^a \|_{\EE^1(T)}  \le  K^a. 
\end{equation}
\end{lem}

The next two propositions are proved in Section~6. 
\begin{prop}
\label{p3} 
There is $C_1 > 0$, such that for all $k \ge k_0$, all $T \le 1$ and all 
$F \in \FF^2(T)$, the Cauchy problem $\eqref{26}$ has a 
unique solution $U \in \EE^2(T)$  and 
\begin{equation}
\label{322}
\|  U  \|_{\EE^2(T)}  \le  C_1 \|  F  \|_{\FF^2(T)} . 
\end{equation}
\end{prop}

The nonlinearity $N_k (U)$ occurring in \eqref{25} is quadratic. Denote by 
$\cN_k  (U, V)$ the bilinear associated  form such that 
$N_k(U) = \cN_k(U, U)$. 

\begin{prop}
\label{p4} 
There is $C_2 > 0$, such that for all $k \ge k_0$, all $T \le 1$ and all 
$U$ and $V$ in  $\EE^1(T)$,   there holds 
   $\cN_k (U, V) \in \FF^2(T)$  and 
\begin{equation}
\label{322bis}
\|  \cN_k (U, V)   \|_{\FF^2(T)}  \le  C_2 \|  U  \|_{\EE^1(T)}   \| V  \|_{\EE^1(T)} . 
\end{equation}
\end{prop}

These estimates easily imply the following: 
\begin{cor}
\label{cor43}
 There are  $c_0 > 0$, $C$  and  $k_0$, such that for all $k \ge k_0$ 
and all $\delta \in ]0, 1]$, the problem $\eqref{28y}$ has a unique solution 
$u = (e, n)$ in the unit ball of  $ \EE^1(T) $, provided that 
\begin{equation}
\label{422}
\delta k^{ -\qt} e^{\sigma T}  \le  c_0. 
\end{equation}

Moreover, the solution satisfies 
\begin{equation}
\label{326}
\| n(t)  \|_{H^s}  \le  C  k^{-\qt}  e^{ \sigma t }
\end{equation}

\end{cor}

\begin{proof}
 Denote by $L_k^{-1} F$ the solution of \eqref{26}, and consider the mapping 
$$
u \ \mapsto \    \cT  u   :=  \delta  L_k^{-1} N_k (u^a + u) 
$$
which, by the lemma and propositions above, is well defined from  $\EE^1(T) $ to $\EE^1(T)$.
Moreover, 
$$
\Vert \cT u \Vert_{\EE^1(T) } \le   
C_1 C_2  \delta k^{-\qt} e^{ \sigma T} (K^a + \Vert u \Vert_{\EE(T)} )^2 . 
$$
Thus it maps the unit ball to of $\EE^1(T)$ to itself, if \eqref{422} holds with $c_0$ small enough. 
Similarly, decreasing $c_0$ if necessary, one shows that this mapping is contractive on the unit ball, 
implying the existence and uniqueness of the solution  of $u = \cT u$ in the unit ball.

The equation $u = \cT u$ and the estimates also imply that 
$$
\begin{aligned}
\| n(t)  \|_{H^s}  \le k^{-\mez} e^{ 2 \sigma t} \| u \|_{\EE^2(T)} 
& \le  C_1 C_2   \delta k^{-\mez }  e^{2 \sigma t }   (K^a + 1 )^2 
\\
 &\le  C_1 C_2   c_0  k^{-\qt  }  e^{ \sigma t }   (K^a + 1 )^2  
\end{aligned}
$$
finishing the proof of the Corollary. 
\end{proof}

We end this section by proving that the main Theorem~\ref{mainth} is a consequence 
of this analysis. 

\medbreak

\begin{proof}[Proof of Theorem~$\ref{mainth}$]

We fix an integer $s$.  
With 
\begin{equation}
\delta  = k^{ - (2s + 2)} , 
\end{equation} 
 Corollary~\ref{cor43} provides us with solutions of \eqref{25}, 
$  U_k =  \underline U  +  \delta (U^a + u_k ) $,   
with $u_k$ in the unit ball of $\EE^1(T_k)$  and 
$T_k = \frac{1}{\sigma} \ln (k^{2s + 2 + \qt} / c_0)$ satisfies 
\begin{equation}
\label{327}
\delta k^{ -1/4} e^{\sigma T_k }  =   c_0.   
\end{equation}
Since $\sigma $ is of order $k^\mez$ by \eqref{310},  $T_k $ tends to 
$0$ as $k$ tends to infinity. 

Going back to the $(z, x)$ variables, according to the change of variables 
\eqref{21}, we obtain solutions,  denoted by $\tilde U_k = \underline U + \tilde u_k$, of the original Zakharov system \eqref{1}. 
 Set  $\tilde u_k  = (\tilde e_k, \tilde n_k)$; these functions are deduced from
 $ \delta (U^a + u_k ) $ by the change of variables \eqref{21}.  
Since $m \le k^2 + k$, we can evaluate the $H^s$ norm (in the variables $(z, x)$) 
of the Cauchy data 
$$
\begin{aligned}
\| (  \tilde e_k{} _{| t = 0}  ,  \tilde n_k{} _{| t = 0} , \D_t  \tilde n_k{} _{| t = 0} \|_{H^s(\TT)} 
&
\le  C \,  \delta  \,   k^{ 2s +1}  \, \| U^a + u_k \| _{\EE^1(T)}
\\
&     \le  
 C \,  \delta  \,   k^{ 2s +1}   (K^a + 1). 
\end{aligned} 
$$
Note that there is no Jacobian factor because the $L^2$ norms are taken for 
$(z, x) \in \TT$ in the left hand side and for $\theta \in \RR/  2\pi \ZZ$ in the right hand side so that
\begin{equation}
\label{328}
\int_\TT v(k x - m z) dzdx = \frac{\mathrm{meas} \TT }{2\pi} \int_{0}^{2 \pi}   v(\theta) d\theta . 
\end{equation} 
 Therefore, with our choice of $\delta$, 
the left hand side tends to zero as $k$ tends to infinity. 

Finally we compute the $L^2$ norm of $\tilde n_k $ at time $T_k$.  Using 
\eqref{328} and   \eqref{37}  we see that 
$$
\| \tilde n_k (T_k) \|_{L^2(\TT)}   \ge    c_1 
\delta  \sinh (T_k  \sigma) \  -  \ \delta  \| n_k (T_k)  \|_{L^2}   
$$
with $c_1 > 0 $ independent of $k$. Therefore,     \eqref{326}  \eqref{327} imply that 
$$
\begin{aligned}
\| \tilde n_k (T_k) \|_{L^2(\TT)}  &  \ge    \mez c_1 
\delta e^{ \sigma T_k}   \  -  \  C \delta   k^{- \qt } e^{ \sigma T_k}    - O (\delta e^{ - \sigma T_k} ) 
\\
& \ge  \mez  c_1 c_0  k^\qt -   C c_0  -  o(1) . 
\end{aligned}
$$
Therefore this $L^2$ norm tends to $+ \infty$ and the proof of the theorem is complete. 
\end{proof}


 \section{The linear instability}

We study  the linear equation for $U = (e, n)$ and $F = (f, g)$
 \begin{equation}
 \label{41a}
 L_k U = F
 \end{equation}
using Fourier series expansions in $\theta$: 
\begin{equation}
\label{27}
e ( \theta, t) = \sum \hat e_p(t) e^{ i p \theta},  \quad n ( \theta, t) = \sum \hat n_p(t) e^{ i p \theta}. 
\end{equation}
Since $n$ and $g$ are   real, 
\begin{equation}
\label{43a} 
\hat n_{-p} = \overline {\hat n_p}, \quad  \hat g_{-p} = \overline {\hat g_p},
\end{equation}    
and \eqref{41a} reduces to 
\begin{equation}
\label{28}
\tilde L_k( \D_t, 0) U_0 := \begin{pmatrix}
      \D_t \hat e_0  - E_0 \hat n_0      \\
      \D^2 _t \hat n_0   
\end{pmatrix}= F_0 := \begin{pmatrix}
      \hat f_0   \\
      \hat g_0  
\end{pmatrix} 
\end{equation}     
and for $p \ge 1$ 
\begin{equation}
\label{29}
\left\{ \begin{aligned}
&    ( i \D_t  +  m p    -  k^2 p^2) \hat  e_p   -  \uE  \hat n_p    =      \hat f_p ,
\\
&  (  i  \D_t +   m p  +    k^2 p^2)  \tilde  e_{p}  +   \overline {\uE} \hat  n_p    =    \tilde f_p ,\\
& (\D_t^2  +  k^2 p^2  )   \hat n_p   + k^2 p^2 \big(    \overline {\uE }  \hat e_ p   
+ \uE \tilde e_p ) 
  =  \hat g_p  , 
\end{aligned}
\right.
\end{equation}
with 
\begin{equation}
\label{210}
\tilde e_p = \overline {e _{-p}}, \quad \tilde f_p = - \overline { f_{-p}}
\end{equation}
are the Fourier coefficients of $\overline e$ and $- \overline f$ respectively. 
For $p > 0$, 
we denote by $\widetilde L_k (\D_t, p)$ the linear operator in the left hand side of 
\eqref{29}. 

\bigbreak
In the remaining part of this section we concentrate on the case $p = 1$ 
and prove Proposition~\ref{p1}.  
We reduce \eqref{29} for $p= 1$  to a first order 
system by  introducing $ v_1  = -  i  k^{-1} \D_t \hat n_1 $. The equation  reads
\begin{equation}
\label{37x}
i \D_t V_1  + A V_1  =   F_1 , 
\end{equation}
 with 
 \begin{equation*}
\label{37ax}
V_1 =  (\hat e_1, \tilde e_1, \hat n_1, v_1) , \quad 
F_1 = (\hat f_1, \tilde f_1, 0, k^{-1} \hat g_1)
\end{equation*}
 and 
 \begin{equation*}
\label{37bx}
A =  \begin{pmatrix}
     m  - k^2      & 0 & -\uE    & 0 \\
     0  &      m  + k^2    & \overline {\uE} & 0 
     \\
     0 & 0 & 0 & k  
     \\
    k    \overline {\uE}   &  k       \uE  & k    & 0 
\end{pmatrix}
\end{equation*}

\begin{lem} 
\label{lem31}
If $\uE \ne 0$ and $k$ is large enough, 
$A$ has four distinct eigenvalues; two, called $\lambda_1$ and $\lambda_2$ are real and the other two, $\lambda_3$ and $\lambda_4$, 
are non-real and complex conjugated. There holds
\begin{equation}
\label{38x}
\lambda_1  \sim 2 k^2, \quad   \lambda_2 \sim - k, \quad  
\re  \lambda_3 \sim   k , \quad \sigma := \im \lambda _3 \sim \vert \uE \vert \sqrt {k/2} . 
\end{equation}

\end{lem}
\begin{proof}
This follows from the analysis of the determinant equation in Section~2. 
The eigenvalue equation is 
\begin{equation*}
\label{4bis}
 P = (   \lambda^2  - k^2) \big(   (\lambda - m)^2 - k^4 \big) - 2 
 \vert \underline E \vert^2 k^4  =  0 
\end{equation*}
Following  \eqref{23}, we write $m = k^2 + k +m'$, 
and the equation reads 
$$
 (   \lambda^2  - k^2)   (\lambda -  k + m')  (\lambda - 2 k^2 - k + m' )  =   2 
 \vert \underline E \vert^2 k^4   
$$
 Because 
 $m' = O(1)$, the lemma easily follows by perturbation analysis of the roots of 
 $$
 (   \lambda^2  - k^2)   (\lambda -  k + m')  (\lambda - 2 k^2 - k + m' )  =   0.
 $$
 \end{proof}

Next,  to evaluate $e^{ it A}$, we need to analyze the eigenprojectors of $A$. 
Denote by $r_j$ [resp. $l_j$] right [resp. left ] eigenvectors of $A$ 
associated to the eigenvalue  $\lambda_j$. Then 
\begin{equation}
\label{36y}
e^{ i t A} \Phi = \sum_{j=1}^4 e^{i t \lambda_j}  \frac{( \l_j \cdot \Phi)}
{(l_j \cdot r_j)}  \,  r_j . 
\end{equation}
A detailed inspection of the eigenvector equations implies the following
\begin{equation*}
r_1 =   \begin{pmatrix} O(k^{-4}) \\ 1 \\ O(k^{-2}) \\ O(k^{-1} ) \end{pmatrix}
, \quad
l_1 =  \begin{pmatrix} O(k^{-4}), & 1,  & O(k^{-2}), &  O(k^{-3} ) \end{pmatrix}, 
\end{equation*}
\begin{equation*}
r_2 \sim  \begin{pmatrix} O(k^{-1}) \\ O(k^{-2} )  \\ 1 \\  -1    \end{pmatrix}
, \quad
l_2 \sim  \begin{pmatrix} O(1), & O(k^{-1}),  & 1 , &  -1   \end{pmatrix}, 
\end{equation*}
where for  vectors, $ a \sim b $ means that all the components 
satisfy  $a_k \sim b_k$. Moreover,  
 \begin{equation*}
r_3 \sim  \begin{pmatrix}  i \uE  / \sigma  \\ O(k^{-2} )  \\ 1 \\  1    \end{pmatrix}
, \quad
l_3 \sim  \begin{pmatrix} k \overline {\uE} / i \sigma , & O(k^{-1}),  & 1 , &  
1   \end{pmatrix}, 
\end{equation*}
 \begin{equation*}
r_4  = \overline r_3 
, \quad
l_4  = \overline l_3 
\end{equation*}
where $\sigma^2 = k \vert \uE  \vert^2/2 \approx k$. 
Note that $r_{3, 4} = 0(1)$ and $r_3 - r_4 = O(|\uE | / \sqrt k)$ and  
$l_{3,4} = O(|E| \sqrt k)$ while 
$r_{3, 4} \cdot l_{3, 4} \sim 4 $. This reflects that  for $\uE = 0$, the corresponding 
matrix has a Jordan block. 
\bigskip

\begin{proof}[Proof of Proposition~$\ref{p1}$.]
{\quad  } 

 With notations as above, 
\begin{equation*}
V_1^a =  \begin{pmatrix}   \hat e^a_1  \\ \tilde e_1^a   \\ \hat n_1^a  \\  v_1^a    \end{pmatrix} :=   \frac{1}{4} \big(   e^{ i t \lambda_4 } r_4   -   e^{ i t \lambda_3 } r_3   \big)  
\end{equation*}
is a solution of 
\eqref{37x} with $F_1 = 0$.  It corresponds to a solution
$(\hat e_1^a, \tilde e_1^a, n_1^a)$ of $\tilde L_1 \widetilde U_1^a  = 0$ and therefore 
to a solution 
  \begin{equation*}
 e^a    = \hat e^a_1   e^{i \theta}  + \overline {\tilde e^a_1}   e^{ - i \theta}, 
\quad  n^a   = \hat n^a_1   e^{i \theta}  + \overline {n^a_1}   e^{ - i \theta}
 \end{equation*}
 of $L_1 U^a = 0$. 
 
 Choosing, as we may,  $r_3$ and $r_4$ such that the third component 
 is exactly equal to one, we obtain that 
 $$
  n^a (t, \theta )  =  \sinh \big(   t \sigma \big)  \cos \big(  t \re \lambda_3  + \theta \big)  
 $$
 and the estimate \eqref{39a} follows from the estimates of the eigenvectors above. 
 Moreover, \eqref{310a} follows from Lemma~\ref{lem31}. 
 \end{proof}

\bigbreak

Next we turn to the analysis of \eqref{37x}. The solution with vanishing initial data
is
\begin{equation}
\label{310z}
V_1(t)  = \sum_{j=1}^4 \int_0^t e^{i (t - s)  \lambda_j}  \frac{( \l_j \cdot F_1(s) )}
{(l_j \cdot r_j)}   \,  r_j  \, ds. 
\end{equation}
Introduce $\Phi_j = l_j \cdot F_1$.  With $f $ denoting $(\hat f_1, \tilde f_1)$
and $g = \hat g_1$
there holds
\begin{equation}
\label{311z}
\begin{aligned}
\Phi_1  &=  * f  + * k^{-4} g, \\
\Phi_2 & =  * f  + * k^{-1} g,  \\
\Phi_{3,4} & =  * \sqrt k f  + * k^{-1} g . 
\end{aligned}
\end{equation} 
where $*$ denotes constants coefficients that are uniformly bounded in $k$. 
Let 
 \begin{equation*}
 \Psi_j(t) = \int_0^t e^{ i \lambda_j (t- s) } \Phi_j (s) ds. 
 \end{equation*}
 The properties of the $r_j$ 's and \eqref{310z} imply that the components 
 $(\hat e_1, \tilde e_1, \hat n_1, v_1)$ of $V_1$ satisfy: 
 \begin{equation}
\label{312z}
\begin{aligned}
\hat e_1  &=  * k^{-4} \Psi_1   + * k^{-1} \Psi_2 +  *  k^{ -1/2} \Psi_{3,4}, \\
\tilde e_1  & =  * \Psi_1    + * k^{-2} \Psi_2  + k^{-2} \Psi_{3, 4},  \\
\hat n_1  & =  * k^{-2} \Psi_1   +  * \Psi_2  +  * \Psi_{3,4} , 
\\ v _1  & =  * k^{-1} \Psi_1   +  * \Psi_2  +  * \Psi_{3,4}   . 
\end{aligned}
\end{equation}

We use the following elementary estimates:
\begin{lem}
\label{lem42}
Let 
 \begin{equation*}
 \psi(t) = \int_0^t e^{ i \lambda (t- s) } \phi_j (s) ds. 
 \end{equation*} 
 There holds
   \begin{equation*}
  \begin{aligned}
\vert  \psi (t) \vert  & \le  \int_0^t e^{ - \im  \lambda  (t- s) } \vert  \phi (s) \vert ds,  
\\
  \vert \D_t \psi (t) \vert  &   \le  \vert \lambda_j \vert \,  \vert  \psi (t) \vert + 
\vert   \phi  (t) \vert, 
 \\
\vert  \D_t\psi (t)&  \vert \le  e^{ - \im \lambda t} \vert   \phi  (0) \vert + 
 \int_0^t e^{ - \im  \lambda  (t- s) } \vert \D_t   \phi  (s) \vert ds,  
 \\
\vert \lambda_j \vert \,  \vert  \psi (t) \vert  &  \le \vert \D_t \psi (t) \vert + 
\vert   \phi  (t) \vert . 
\end{aligned}
 \end{equation*}
\end{lem}

To simplify notations, we note 
$A \ls B $ to mean that  there is a constant $C$ independent of $k$ such that 
$A \le C B $. 
  We use the first and second estimate of Lemma~\ref{lem42} to bound the contributions of 
  $g$ to the integrals in \eqref{310z}, and we use
  the third and fourth estimate, when necessary, to bound the contributions of $f$. 
  Therefore, 
  \begin{equation*}
  \begin{aligned}
  \vert \Psi_1(t) \vert  & \ls \int_0^t \vert f(s), k^{-4} g(s) \vert  ds  , 
  \\
   \vert \D_t \Psi_1(t) \vert & \ls  \vert f(0) \vert + \vert k^{-4}  g(t) \vert + \int_0^t \vert \D_t f(s), k^{-2} g(s) \vert  ds  , 
   \\
  k^2  \vert    \Psi_1(t) \vert & \ls  \vert f(0) \vert + \vert   f(t) \vert + \int_0^t \vert \D_t f(s), k^{-2} g(s) \vert  ds 
  \end{aligned}
  \end{equation*}
   \begin{equation*}
  \begin{aligned}
  \vert \Psi_2(t) \vert  & \ls \int_0^t \vert f(s), k^{-1} g(s) \vert  ds  , 
  \\
   \vert \D_t \Psi_2(t) \vert & \ls  \vert f(0) \vert + \vert k^{-1}  g(t) \vert + \int_0^t \vert \D_t f(s),  g(s) \vert  ds  ,
   \\
  k   \vert  \Psi_2(t) \vert & \ls  \vert f(0) \vert + \vert f(t) \vert + \int_0^t \vert \D_t f(s),  g(s) \vert  ds  , 
  \end{aligned}
  \end{equation*}
   \begin{equation*}
  \begin{aligned}
  \vert \Psi_{3,4}(t) \vert  & \ls \int_0^t 
  e^{ (t-s)\sigma } \vert \sqrt k f(s), k^{-1} g(s) \vert  ds  , 
  \\
   \vert \D_t \Psi_{3,4}(t) \vert & \ls e^{t \sigma }  \vert \sqrt k f(0) \vert + \vert k^{-1}  g(t) \vert + \int_0^t   e^{ (t-s)\sigma }  \vert \sqrt k  \D_t f(s),  g(s) \vert  ds, 
   \\
   k   \vert  \Psi_{3,4}(t) \vert & \ls e^{t\sigma}  \vert \sqrt k f(0) \vert +
    \vert \sqrt k  f (t) \vert + \int_0^t   e^{ (t-s)\sigma } 
     \vert \sqrt k  \D_t f(s),  g(s) \vert  ds .  
  \end{aligned}
  \end{equation*}

 Adding up the various estimates, we obtain: 
 
 \begin{prop} 
 \label{prop34}
For $p = 1$, the solution 
  $ (\hat e_1, \tilde e_1, \hat n_1)$ of $\eqref{29}$    with vanishing 
  initial data  satisfies: 
  \begin{equation}
  \label{314z}
  \begin{aligned}
 & \vert \hat e_1 (t) \vert \ls \int_0^t e^{ \sigma (t-s) } \vert   f_1(s) , k^{-\frac{3}{2}} \hat g_1(s) \vert ds, \\ 
  & \begin{aligned}
  \vert \D_t \hat e_1(t) \vert \lesssim    
  e^{\sigma t  } \vert f_1 (0) \vert &+    \vert k^{-\frac{3}{2}} \hat g_1 (t) \vert 
  \\ & +      \int_0^t  e^{\sigma (t-s)   } \vert   \D_t f_1(s), k^{-\mez} \hat g_1(s)  \vert  ds, 
  \end{aligned}
\end{aligned}
  \end{equation} 
  \begin{equation}
  \label{315z}
  \begin{aligned}
k^2   \vert \tilde e_1 (t) \vert + \vert \D_t \tilde e_1(t) \vert \lesssim  & 
  e^{ \sigma  t  } \vert f_1 (0) \vert  +  \vert f_1(t) \vert +   \vert k^{-3} \hat g_1 (t) \vert  \\
  &  +      \int_0^t  e^{ \sigma  (t-s)  } \vert  \D_t f_1(s), k^{-1} \hat g_1(s)  \vert  ds,  
\end{aligned}
  \end{equation} 
   \begin{equation}
   \label{316z}
  \begin{aligned}
  k   \vert \hat n_1 (t) \vert + \vert \D_t \hat n_1(t) \vert \lesssim  & 
  e^{ \sigma t  }   k^\mez  \vert f_1 (0) \vert  +    \vert k^\mez  f_1 (t) \vert +   \vert k^{-1}  \hat g_1 (t) \vert  \\
  &  +      \int_0^t  e^{ \sigma  (t-s)   } \vert f_1  (s),  k^\mez \D_t f_1(s), \hat g_1(s)  \vert  ds,  
\end{aligned}
  \end{equation} 
  where $f_1 = (\hat f_1, \tilde f_1)$. 
 \end{prop}
 
 \begin{cor}
 \label{cor44} 
 There are $k_0$ and $C$ such that for all $k \ge k_0$,  $K$, 
 $T > 0$, and all 
 $ f_1 = (\hat f_1, \tilde f_1) $, $ g_1$ satisfying for $t \in [0, T]$ 
$$
    k^\mez  |   f_1 (t) |  +   k^{- \mez}     |  \D_t  f _1  (t) |  + k^{-\tqt} | \hat g_1(t) |  
  \le K  e^{2 \sigma t }  ,
$$
then the solution of $\eqref{29}$ for $p= 1$  with vanishing initial data satisfies 
\begin{eqnarray*}
& k   \vert \hat e_1  (t) \vert  +       \vert \D_t \hat e_1  (t) \vert  
  \le   C K   e^{2 \sigma t } ,
\\
&   k  | \tilde  e_{-1}  (t) |  +    k^{- \qt }  |  \D_t  \tilde e_{-1}(t) | 
  \le C K  e^{ 2 \sigma t } ,
\\
&
  k^\mez   |  \hat n_1  (t) |  +  k^{- \mez}  |  \D_t \hat n_1  (t) | 
  \le C  K e^{2 \sigma t } . 
\end{eqnarray*}
 \end{cor}

\begin{proof}
{\bf a) } From Proposition~\ref{prop34} we deduce that 
\begin{equation*}
k  \vert \hat e_1 (t) \vert \le C  K_1 \sqrt k 
 \int_0^t e^{ \sigma (t-t') } e^{ 2 \sigma t'}  dt' \le C  K  e^{ 2 \sigma t} , 
\end{equation*}
where we have used that $\sigma \approx \sqrt k$. Similarly, 
\begin{equation*}
  \vert \D_t \hat e_1(t) \vert \le C K_1   \Big( 
 k^{-1/2}  e^{\sigma t  }  +     k^{-3/4} e^{ 2 \sigma t  } 
 +        \int_0^t  \sqrt k e^{ \sigma (t-t') } e^{ 2 \sigma t'}  dt' \Big) \le C  K  e^{ 2 \sigma t}. 
  \end{equation*} 
  This implies the first estimate. 
  
  \medbreak
  {\bf b) } Similarly, \eqref{315z} implies that 
    \begin{equation}
    \label{414}
  \begin{aligned}
k^2   \vert \hat  e_{-1} (t) \vert + \vert \D_t \hat e_{-1}(t) \vert \le C K_1  \Big( 
  e^{\sigma t  }  +     e^{ 2 \sigma t  } 
 +        \int_0^t  \sqrt k e^{ \sigma (t-t') } e^{ 2 \sigma t'}  dt' \Big)\\
  \le C K e^{ 2 \sigma t}.  
\end{aligned}
  \end{equation} 

{\bf c) } The estimate \eqref{316z} implies that 
     \begin{equation}
     \label{417}
  \begin{aligned}
  k   \vert \hat n_1 (t) \vert + \vert \D_t \hat n_1(t) \vert \le     C K_1  \Big( 
   e^{ \sigma t} +   e^{ 2 \sigma t  }  +    \    \int_0^t    k e^{ \sigma (t-t') } e^{ 2 \sigma t'}  dt' \Big)\\
  \le C  K \sqrt k e^{ 2 \sigma t} 
\end{aligned}
  \end{equation} 
  and the lemma is proved.
\end{proof}

 \section{The linear equation}

 We continue the analysis of the linear equation \eqref{41a}.  
 As seen in \eqref{29}, when expanded in Fourier series, this equation couples
 the coefficients of indices $p$ and $-p$.  The case of indices $+1$ and $-1$
 is studied in the previous section. 
 Using the notations 
 \begin{equation}
 v = \hat v_1 e^{ i \theta} + \hat v_{-1} e^{ - i \theta}  + v'' 
 \end{equation}
 we consider the equation \eqref{41a} for functions with vanishing Fourier coefficients 
 of indices $\pm 1$ : 
 \begin{equation}
 \label{51a}
 L''_k U'' = F'' , 
 \end{equation}
which reduces to the analysis of  equations \eqref{29} for Fourier   $p  \ne 1$.   
  
The symbol of $\widetilde L_k(\D_t, p) $ is 
\begin{equation}
\label{211}
\widetilde L_k( i \tau, p) =  \begin{pmatrix}
   -  \tau+ mp - k^2 p^2    & 0 & -E_0   \\
     0  &  -   \tau  + mp + k^2 p ^2 & \overline E_0  
     \\
    k^2 p^2  \overline E_0  &  k^2 p^2    E_0  & k^2 p^2  - \tau^2 
\end{pmatrix}\end{equation}
which is of course equal to  the symbol \eqref{3} with 
 with $\xi = k p$,  $\zeta = - m p $,  up to a change of sign in the first line.

\medbreak

Assume first that  $p > 1$. 
 In this case, we consider $\widetilde L_k(\D_t, p)$ as a perturbation of 
 \begin{equation}
\label{212}
M_k(\D_t, p) := \begin{pmatrix}
        i \D_t  +  m p    -  k^2 p^2     \\
       i  \D_t +   m p  +    k^2 p^2
       \\  \D_t^2  +  k^2 p^2 
\end{pmatrix}
\end{equation} 
  For the wave operator, we use the classical estimates:  
 
 \begin{lem} There is $C > 0$, such that for all $k \ge 1$ and 
 $p \ge 1$, the solution $n$ of 
 \begin{equation}
\label{38}
\D_t^2 n + k^2 p^2 n = g, \quad n(0) = \D_t n (0) = 0 
\end{equation}
 satisfies
 \begin{equation}
\label{39}
 k p \vert  n(t) \vert  + | \D_t n(t) | \le C \Vert g \Vert_{L^1([0, t])}. 
\end{equation}
 \end{lem}
 
 For the Schr\"odinger equations, we use the following estimates. 
 
  \begin{lem} There are $C > 0$ and $k_0 \ge 1$, such that for all $k \ge k_0$ and 
 $p \ge 2$, the solutions  of 
 \begin{equation}
\label{310}
(  i \D_t + mp  \pm  k^2 p^2 )e  = f , \quad e(0)  
\end{equation}
 satisfy
 \begin{equation}
\label{311}
 k^2 p^2 \vert  e(t) \vert  + | \D_t e(t) | \le C\big(   \Vert   f  \Vert_{L^1([0, t])} + 
 \Vert \D_t f  \Vert_{L^1([0, t])} + | f(0) | \big)  
\end{equation}
 \end{lem}
\begin{proof}
Standard energy estimates imply that 
\begin{equation*}
\label{311a}
 \vert  e(t) \vert  \le  C  \big( |e(0)| +    \Vert f  \Vert_{L^1([0, t])} \big) . 
\end{equation*}
Differentiating in time the equation, we obtain 
\begin{equation*}
\label{311b}
 \vert  \D_t e(t) \vert  \le  C  \big( |\D_t e(0)| +    \Vert \D_t f  \Vert_{L^1([0, t])} \big) . 
\end{equation*}
The initial condition  in \eqref{310} implies that $\D_t e (0) = -i f(0)$. 
Therefore, 
\begin{equation*}
\label{311c}
  \vert  ( k^2 p^2 \pm  m p ) e(t) \vert  + | \D_t e(t) | \le C\big(   \Vert   f  \Vert_{L^1([0, t])} + 
 \Vert \D_t f  \Vert_{L^1([0, t])} + | f(0) |  + | f(t) | \big)  
\end{equation*}
Recall that $m$ is linked to $k$ through \eqref{23}. Thus $m p \le k^2 p + k p$ and 
$k^2 p^2 - mp \ge k^2 (p^2 - p ) - k p \ge c k^2 p^2$ for 
all $p \ge 2$ if $k$ is large enough. 
\end{proof}

 \begin{prop}
 \label{prop37}
 Consider the equation $\eqref{29}$ with initial data 
 \begin{equation}
\label{312}
\hat e_p(0) = \tilde e_p(0) = \hat n_p(0) = \D_t \hat n_p(0) = 0
\end{equation}
 Then, for $p \ge 2$, $k \ge k_0$, there holds for $t \in [0,1]$: 
 \begin{equation}
\label{313}
\begin{aligned}
k^2 p^2   | \hat e_p(t) ,  \tilde e_p(t) |   +     |\D_t  & \hat e_p(t) , \D_t  \tilde e_p(t)   |  + 
kp | \hat n_p(t) | + | \D_t \hat n_p(t) | 
\\    \le &  C \big(  \Vert   \hat f_p, \tilde f_p    \Vert_{L^1([0, t])}    +  
 \Vert \D_t  \hat f_p, \D_t \tilde f_p   \Vert_{L^1([0, t])} 
 \\ &  + |  \hat f_p(0), \tilde f_p (0) |  + |  \hat f_p(t), \tilde f_p (t) | + \Vert   \hat g_p    \Vert_{L^1([0, t])} 
 \big). 
 \end{aligned} 
\end{equation}
 \end{prop} 
 \begin{proof}
 The lemmas above imply that the left hand side is estimated by the right hand side plus 
 \begin{equation*}
\label{313a}
C\Big(   | \hat n_p(t)| + \Vert \hat n_p, \ \D_t \hat n_p , \ k^2 p^2     \hat e_p(t) , \ k^2 p^2   \tilde e_p  
\Vert_{L^1([0, t])}\Big)
\end{equation*}
The first term is absorbed in the left hand side by $k p  |\hat n_p(t)| $ for 
$k$ large enough. With Gronwall's lemma, this implies \eqref{313} for $t \in [0, 1]$, with a larger constant $C$. 
 \end{proof}

 When $p = 0$, there holds: 
 \begin{lem}
 \label{lem38}
 When $\hat g_0 = 0$, the solution of $\eqref{28}$ with vanishing initial data is
 \begin{equation}
\label{324}
\hat n_0 = 0,  \quad  \hat e_0 (t) = \int_0^t \hat f_0(t') dt'. 
\end{equation}
 \end{lem} 

With the estimates \eqref{313}, one deduces the following result 

\begin{cor}
\label{cor55}
 There are $k_0$ and $C$ such that for all $k \ge k_0$,  $K$, 
 $T > 0$, and all 
 $ (f'' , g'') $ with $\hat g_0 = 0$,  satisfying for $t \in [0, T]$ 
\begin{eqnarray*}
& k^\mez   \Vert f'' (t) \Vert_{H^s} +   k^{- \mez}     \Vert \D_t  f '' (t) \Vert_{H^s} 
  \le K   e^{2 \sigma t } ,
\\
&     \Vert g ''  (t) \Vert_{H^s}   
  \le K    k^{3/4}    e^{2 \sigma t } . 
\end{eqnarray*}
the solution of $\eqref{51a}$ with vanishing initial data  satisfies 
\begin{eqnarray*}
&   k  \Vert e'' (t) \Vert_{H^{s+2}} +    k^{-\qt}   \Vert \D_t  e'' (t) \Vert_{H^s}
  \le C  K  e^{ 2 \sigma t } ,
\\
&
  k^\mez  \Vert n'' (t) \Vert_{H^{s }} + k^{- \mez}      \Vert \D_t n'' (t) \Vert_{H^s} 
  \le C  K  e^{2 \sigma t } . 
\end{eqnarray*}
\end{cor} 

\begin{proof}

  By Lemma~\ref{lem38}, there holds
     \begin{equation}
    \label{415}
k  \vert \hat e_0 (t) \vert + \vert \D_t \hat e_0(t) \vert \le C K_1  \Big( 
     e^{ 2 \sigma t  } 
 +        \int_0^t  \sqrt k e^{ \sigma (t-t') } e^{ 2 \sigma t'}  dt' \Big) 
  \le C  K e^{ 2 \sigma t}.  
  \end{equation} 
  Next, Proposition~\ref{prop37} implies that 
 $e''$ 
  satisfies 
   \begin{equation}
    \label{416}
  \begin{aligned}
k^2    & \Vert \D_\theta^2  e'' (t) \Vert_{H^s}  + \Vert \D_t   e'' (t) \Vert_{H^s}
\\ &  \le C K   \Big( 
    ( 1  +     e^{ 2 \sigma t  } ) 
 +        \int_0^t   k^{3/4}  e^{ \sigma (t-t') } e^{ 2 \sigma t'}  dt' \Big) 
  \le k^{1/4}  C  K  e^{ 2 \sigma t}.  
\end{aligned}
  \end{equation} 
Together with \eqref{415}  this implies  the first estimate. 

Moreover, Proposition~\ref{prop37} implies that  $ n''$ 
satisfies
\begin{equation}
    \label{418}
  \begin{aligned}
k    & \Vert n'' (t) \Vert_{H^s}  + \Vert \D_t   n'' (t) \Vert_{H^s}
\\ &  \le C K   \Big( 
    ( 1  +     e^{ 2 \sigma t  } ) 
 +        \int_0^t   k^{3/4}  e^{ \sigma (t-t') } e^{ 2 \sigma t'}  dt' \Big) 
  \le k^{1/4}  C K  e^{ 2 \sigma t}.  
\end{aligned}
  \end{equation} 
  Since $\hat n_0 = 0$, this  implies the second estimate . . 

\end{proof}


\section{End of proofs}

   First, we note that  Proposition~\ref{p3} is an immediate consequence of Corollaries~\ref{cor44} and 
  \ref{cor55}. 
  
  It remains to prove Proposition~\ref{p4}.  
  With $U = (e, n)$ and $  U^*    =  ( e^* ,  n^*)$, there holds 
  \begin{equation}
  \label{61}
  \cN_k (U,   U^* )  = (f, g) with 
  \end{equation}
\begin{eqnarray}
\label{62}
&& f =    n     e^*  +    n^*   e , 
\\
&&   g =   k^2 \D_\theta^2 \big\{   \re (  \overline e   e^*) \big\}. 
\end{eqnarray}
Proposition~\ref{p4} follows from the next estimates.

\begin{lem}
\label{lem41}
There is a constant $C $, independent of $k$, such that 
\begin{eqnarray}
\label{46}
&&\sqrt k  \Vert f (t) \Vert_{H^s} +  \frac{1}{\sqrt k}     \Vert \D_t  f  (t) \Vert_{H^s} 
  \le C   e^{2 \sigma t } \| U \|_{\EE^1(T)}    \|   U^* \|_{\EE^1(T)}   ,
\\
&&\label{47}     \Vert g  (t) \Vert_{H^s}   
  \le C    k^{3/4}    e^{2 \sigma t }  \| U \|_{\EE^1(T)}    \|   U^* \|_{\EE^1(T)}  . 
\end{eqnarray}
Moreover, the mean value $\hat g_0$ of $g$ vanishes. 
\end{lem} 

\begin{proof}
The first estimate follows directly from the definitions  and 
the inequality
\begin{equation*}
\Vert a b \Vert_{H^s} \le C \Vert a \Vert_{H^s} \Vert b \Vert_{H^s}. 
\end{equation*}
Next, we note that for $e = \hat e_1 e^{ i \theta} + e{}' $
and $  e^* =     \hat e^*_1e^{i\theta}  +  e^*{}' $  
\begin{equation*}
 \D_{\theta}^2 (   \overline e  e^* ) 
    =   \D_{\theta}^2 ( \overline e'  e^*{}')   +   
   \overline {\hat e_1} \D_{\theta}^2 ( e^*{} '   e^{ -  i \theta} ) 
   + 
  {\hat e^*_1} \D_{\theta}^2 ( \overline e'   e^{   i \theta} ) . 
\end{equation*}
Hence, in $H^s$ norms, there holds
\begin{equation*}
\begin{aligned}
\Vert   \D_\theta^2  (   \overline e  e^* )    \Vert_{H^s}  
\ls  &    \Vert  \D_\theta^2  e' \Vert_{H^s} \ 
 \big(  \Vert     e^*{}'   \Vert + \Vert \D_\theta e^*{}' \Vert^2 \big)+  
  \Vert  \D_\theta^2  e^*{}' \Vert_{H^s} \ 
 \big(  \Vert     e{}'   \Vert + \Vert \D_\theta e{}' \Vert^2 \big)
  \\
   + &  \vert  \hat e_1 \vert  \big( \Vert \D_\theta^2 e^*{}' \Vert + \Vert  e^*{}' \Vert  \big)  
   +  \vert  \hat e^*_1 \vert  \big( \Vert \D_\theta^2 e{}' \Vert + \Vert  e{}' \Vert  \big)
\end{aligned}
\end{equation*}
 and \eqref{47} follows. 

In addition, the $\theta$-mean value $\hat g_0$ vanishes since $g$ is a $ \theta$-derivative.   
\end{proof}


\vfill \eject

\end{document}